\documentclass{amsart}

\usepackage{amssymb}
\usepackage{amsthm}
\usepackage{epsfig}
\usepackage{color}
\usepackage{enumitem}

\usepackage{latexsym,subfigure}
\usepackage{comment}
\usepackage{amsmath,amsthm,amsfonts,amssymb,graphicx,epsfig,latexsym,color}
\usepackage[pdftex,colorlinks]{hyperref} 

\begin{document}

\newtheorem{theorem}{Theorem}[section]
\newtheorem{lemma}[theorem]{Lemma}
\newtheorem{corollary}[theorem]{Corollary}
\newtheorem{proposition}[theorem]{Proposition}
\newtheorem{question}[theorem]{Question}
\newtheorem{conjecture}[theorem]{Conjecture}
\newtheorem{case}{Case}
\newtheorem{subcase}{Case}[case]
\newtheorem{subsubcase}{Case}[subcase]
\newtheorem{claim}[theorem]{Claim}
\newtheorem{step}{Step}

\theoremstyle{definition}
\newtheorem{remark}[theorem]{Remark}
\newtheorem{definition}[theorem]{Definition}
\newtheorem{example}[theorem]{Example}

\newcommand\bY{\mathbb{Y}}

\newcommand\eps{\epsilon}
\newcommand\E{\mathbb{E}}
\newcommand\Z{\mathbb{Z}}
\newcommand\R{\mathbb{R}}
\newcommand\T{\mathbb{T}}
\newcommand\C{\mathbb{C}}
\newcommand\CC{\mathcal{C}}
\newcommand\N{\mathbb{N}}
\newcommand\G{\mathbf{G}}
\newcommand\A{\mathbb{A}}
\newcommand\HH{\mathbb{H}}
\newcommand\SL{\operatorname{SL}}
\newcommand\Upp{\operatorname{Upp}}
\newcommand\Cay{\operatorname{Cay}}
\newcommand\im{\operatorname{im}}
\newcommand\bes{\operatorname{bes}}
\newcommand\sss{\operatorname{ss}}
\newcommand\GL{\operatorname{GL}}
\newcommand\PGL{\operatorname{PGL}}
\newcommand\SO{\operatorname{SO}}
\newcommand\SU{\operatorname{SU}}
\newcommand\PU{\operatorname{PU}}
\newcommand\PSU{\operatorname{PSU}}
\newcommand\PSL{\operatorname{PSL}}
\newcommand\Rad{\operatorname{Rad}}
\newcommand\Mat{\operatorname{Mat}}
\newcommand\ess{\operatorname{ess}}
\newcommand\rk{\operatorname{rk}}
\newcommand\Alg{\operatorname{Alg}}
\newcommand\Supp{\operatorname{Supp}}
\newcommand\sml{\operatorname{sml}}
\newcommand\lrg{\operatorname{lrg}}
\newcommand\tr{\operatorname{tr}}
\newcommand\Hom{\operatorname{Hom}}
\newcommand\Lie{\operatorname{Lie}}
\newcommand\Inn{\operatorname{Inn}}
\newcommand\Out{\operatorname{Out}}
\newcommand\Aut{\operatorname{Aut}}
\newcommand\diam{\operatorname{diam}}

\renewcommand\P{\mathbb{P}}
\newcommand\F{\mathbb{F}}
\newcommand\Q{\mathbb{Q}}
\renewcommand\b{{\bf b}}
\def\g{\mathfrak{g}}
\def\h{\mathfrak{h}}
\def\n{\mathfrak{n}}
\def\a{\mathfrak{a}}
\def\p{\mathfrak{p}}
\def\q{\mathfrak{q}}
\def\b{\mathfrak{b}}
\def\r{\mathfrak{r}}

\def\bn{\mathbf{n}}

\newcommand{\kf}[1]{\marginpar{\tiny #1 }}
\newcommand{\red}{\textcolor{red}}
\newcommand{\blue}{\textcolor{blue}}

\title{Expansion in simple groups}
\date{\today}

\author{Emmanuel Breuillard}
\address{DPMMS, University of Cambridge, UK}
\email{emmanuel.breuillard@maths.cam.ac.uk}

\author{Alexander Lubotzky}
\address{Department of Mathematics, Hebrew University}
\email{alex.lubotzky@mail.huji.ac.il}
\thanks{ }

\begin{abstract} Two short seminal papers of Margulis used Kazhdan's property $(T)$ to give, on the one hand, explicit constructions of expander graphs, and to prove, on the other hand, the uniqueness of some invariant means on compact simple Lie groups. These papers opened a rich line of research on expansion and spectral gap phenomena in finite and compact simple groups. In this paper we survey the history of this area and point out a number of problems which are still open. 
\end{abstract}

\maketitle

\par\vspace*{.01\textheight}{\centering \emph{Dedicated to Grisha Margulis with admiration and affection} \par}

\setcounter{tocdepth}{1}

\tableofcontents

\section{Introduction}

Grisha Margulis has the Midas touch: whatever he touches becomes gold. It seems that he did not have a particular interest in combinatorics, but in the early seventies events of life brought him to work at the Institute for Information Transmission in Moscow, where he became aware of the concept of expander graphs. Such graphs were known to exist at the time only by counting considerations (\`a la Erd\H{o}s random graph theory), but because of their importance in computer science, explicit constructions were very desirable. Margulis noticed that such explicit constructions could be made using the (new at the time) Kazhdan property (T) from representation theory of semisimple Lie groups and their discrete subgroups. His short paper \cite{margulis73} opened a new area of research with a wealth of remarkable achievements. 

A similar story happened with Margulis' contribution to the so-called Ruziewicz problem, namely: must every rotation invariant finitely additive measure on the sphere $S^n$ be equal to the Lebesgue measure?  It had been known for a long time that the answer is ``no'' for $n=1$, but Margulis \cite{margulis80} (as well as Sullivan \cite{sullivan81}) showed, again using Kazhdan's property (T), that the answer is ``yes'' for $n \geq 4$.

These two seemingly unrelated topics are actually very much connected. This was explained in detail in \cite{lub94}. We will repeat it in a nutshell in \S 2, and give a brief historical description. Both directions of research led to a problem of the following type: 

\bigskip

\noindent {\bf Problem.} \emph{Let $G$ be a non-abelian finite (resp. compact Lie) group. For a finite symmetric, i.e. $S=S^{-1}$,  subset $S$ of $G$, consider $$\Delta^S=\sum_{s \in S} s$$ as an operator on $\mathbf{L}^2(G)$, where $sf(x):=f(s^{-1}x)$. It is easy to see that its largest eigenvalue is $k=|S|$ (with the constant functions being the eigenspace). It has multiplicity one if and only if $S$ generates (resp. generates topologically) $G$. Find $S$ with spectral gap, i.e. for which the second largest eigenvalue of $\Delta^S$ is bounded away from $|S|$. }

This problem has many variants. Do we take $S$ optimal (``best case scenario''), worst (``worst case scenario''), or random? Do we want $k=|S|$ to be fixed? Is the ``bounded away'' uniform? in what? the generators? all groups? all generators?

A quite rich theory has been developed around these questions, which grew out from the above two papers of Margulis. The goal of this paper is to describe this story and to point out several problems that are still open. 

In \S 2, we will give some more history and show how the central problem we study here is related to expanders and to the Ruziewicz problem. 

In \S 3, we describe the numerous developments the expansion problem for finite simple groups has had in the last decade or so. This led to amazing connections with additive combinatorics, diophantine geometry, Hilbert's 5th problem and more. It also led to some new ``non-commutative sieve method'' with some remarkable applications. These subjects have been discussed in a number of books and surveys \cite{lub94,tao-expansion, lubotzky-bulletin-ams,breuillard-minneapolis, breuillard-standrews, breuillard-icm} so we do not cover them here.

In \S 4, we will turn our attention to the compact simple Lie groups. Here much less is known. This direction has recently  received renewed interest from questions in quantum computation (``golden gates''). Now, one looks not only for topological generators $S:=\{g_1,\ldots,g_k\}$ in $G$ with spectral gap but we also want an algorithm that will enable us to find for every $g \in G$ a ``short'' word $w$ in $g_1,.\ldots,g_k$ such that $w(g_1,\ldots,g_k)$ is very close to $g$.

This paper, which only illustrates a small part of Margulis' influence on modern mathematics is dedicated to Grisha with admiration and affection. He has been a personal and professional inspiration for both of us.

\section{Expanders and invariant means}

A family of finite $k$-regular graphs $X_i=(V_i,E_i)$ is called \emph{an expanding family}, if there exists $\epsilon>0$ such that for every $i$ and every subset $Y \subset V_i$ with $|Y| \leq |V_i|/2$, 
$$|\partial Y| \geq \epsilon |Y|,$$
where $\partial Y = E(Y, \overline{Y})$ is the set of edges going out from $Y$ to its complement $\overline{Y}$. 

Margulis made the following seminal observation, which connected expanders and representation theory:

\begin{proposition}[Margulis]\label{margulisprop} Let $\Gamma$ be a group generated by a finite set $S$ with $S=S^{-1}$ and $|S|=k$. Assume that $\Gamma$ has Kazhdan property $(T)$, then the family of finite $k$-regular Cayley graphs $Cay(\Gamma/N; S)$ when $N$ runs over the finite index normal subgroups of $\Gamma$, forms an expanding family.
\end{proposition}

Let us give a sketch of proof: property $(T)$ means that the trivial representation is an isolated point in the unitary dual of $\Gamma$, the space of irreducible unitary representations of $\Gamma$ up to equivalence endowed with the Fell topology. In concrete terms for $\Gamma=\langle S \rangle$ as above, it says that there exists an $\epsilon'>0$, such that whenever $\Gamma$ acts unitarily  on a Hilbert space $\mathcal{H}$ via a (not necessarily irreducible) unitary representation $\rho$ without non-zero fixed vector, for every vector $v \neq 0$ in $\mathcal{H}$, there exists $s\in S$ such that 

\begin{equation}\label{rep} \|\rho(s)v-v\| \geq \epsilon' \|v\|.
\end{equation}

In our situation, let $Y$ be a subset of $\Gamma/N$, i.e. a subset of vertices in $Cay(\Gamma/N,S)$. Let $f$ be the function $f$ in $L^2(\Gamma/N)$ defined by 
$$f(y) = |\overline{Y}|$$ if $y\in Y$ and 
$$f(y)=-|Y|$$  if $y\in \overline{Y}$, where $\overline{Y}$ is the complement of $Y$.

Then $f \in L^2_0(\Gamma/N)$, i.e. $\sum_{y \in V} f(y) = 0$. Now, $\Gamma$ acts unitarily by left translations on $L^2_0(\Gamma/N)$, which, as a representation, is a direct sum of non-trivial irreducible representations. We may thus apply $(\ref{rep})$ and deduce that there exists $s \in S$ such that

\begin{equation}\label{reprep} \|\rho(s)f-f\| \geq \epsilon' \|f\|.
\end{equation}

Spelling out the meaning of $f$, and noting that $f$ is essentially the (normalized) characteristic function of $Y$, one sees that 

$$|sY \triangle Y | \geq \epsilon'' |Y|$$
which implies the desired result. See \cite{lub94} for the full argument with the constants involved.

This fundamental argument gave a lot of families of finite (simple) groups which are expanding families. Every ``mother group'' $\Gamma$ with property $(T)$ gives rise to a family of expanders. 

For example, for every $n \geq 3$, $\Gamma_n = \SL_n(\Z)$ has $(T)$ by Kazhdan's theorem \cite{kazhdan}. Fix a finite set $S$ of generators in $\Gamma$. One deduces that the family $\{Cay(\PSL_n(p),S) ; p \textnormal{ prime}\}$ is a family of expanders. Naturally this raises the question whether all $\PSL_n(p)$ (all $n$ all $p$) or even all non-abelian finite simple groups can be made into an expanding family simultaneously. This will be discussed in \S 3. Meanwhile, let us give it another interpretation. The well-known result of Alon, Milman and others (see \cite{lub94} for detailed history) gives the connection between expanders and spectral gap. Let us formulate it in the context of finite groups.

\begin{proposition} Let $\{G_i\}_{i \in I}$ be a family of finite groups, with symmetric generating sets $S_i$ with $|S_i|=k$ for every $i$. The following are equivalent:

\begin{enumerate}
\item $\{Cay(G_i;S_i)\}_{i \in I}$ forms an expanding family of $k$-regular graphs.
\item There exists $\epsilon'>0$ such that all eigenvalues of 
$$\Delta^{S_i} =  \sum_{s \in S_i} s$$
acting on $L^2_0(G_i)$ are at most\footnote{This is a one-sided condition, i.e. for each eigenvalue $\lambda \leq k-\epsilon'$. However it is also equivalent to the two-sided condition $|\lambda| \leq k-\epsilon'$, provided the $G_i$'s do not have an index $2$ subgroup disjoint from $S_i$, see \cite[Appendix E.]{bggt2} and \cite{biswas}.} $k-\epsilon'$.
\end{enumerate}

\end{proposition}

Let us now move to the invariant mean problem (a.k.a. the Ruziewicz problem) and we will see that the same spectral gap issue comes up. The problem to start with was formulated for the spheres $S^n$ upon which $G=SO(n+1)$ acts. But it generalizes naturally to the group $G$ itself and actually to every compact group, so we will formulate it in this generality.

Let $G$ be a compact group. An invariant mean from $L^\infty(G)$ to $\R$ is a linear functional $m$ satisfying for every $f \in L^\infty(G)$, 

\begin{itemize}
\item $m(f) \geq 0$ if $f \geq 0$,
\item $m(1_G)=1$ and
\item $m(g.f)=m(f)$ for all $g \in G$,
\end{itemize}
where $1_G$ is the constant function equal to $1$ on $G$, and for $g \in G$, $g.f(x)=f(g^{-1}x)$. Integration against the Haar measure, or Haar integration, is such an invariant mean. It is the only such if we assume in addition that $m$ is countably additive. The question is whether this is still true also among the finitely additive invariant means.

\begin{theorem}[Rosenblatt \cite{rosenblatt-transactions} ]\label{rosen} Let $G$ be a compact group, $S=S^{-1}$ a finite symmetric set in $G$ with $|S|=k$ and $\Gamma=\langle S \rangle$ the subgroup generated by $S$. The following are equivalent:
\begin{itemize}
\item The Haar integration is the unique $\Gamma$-invariant mean on $L^\infty(G)$,
\item There exists $\epsilon'>0$ such that all the eigenvalues of $\Delta^S =  \sum_{s \in S} s$ acting on $L_0^2(G)$ are at most $k-\epsilon'$.
\end{itemize}
\end{theorem}

It is not surprising now that if $\Gamma=\langle S \rangle \subset G$ is a dense subgroup with property $(T)$ (i.e. $\Gamma$ has $(T)$ as an abstract discrete group) then every $\Gamma$-invariant mean (and hence $G$-invariant mean) of $G$ is equal to the Haar integration. This was the way Margulis \cite{margulis80} and Sullivan \cite{sullivan81} solved the Ruziewicz problem for $S^n$ (and $SO(n+1)$), $n \geq 4$, to start with. 

Note that when $G$ is finite, then the spectral gap property is not so interesting for a single group $G$ ; it just says that $S$ generates $G$. For an infinite compact group, $S$ generates $G$ topologically if and only if $k$ has multiplicity $1$, i.e. all other eigenvalues are less than $k$. But we want them (there are infinitely many of them!) to be bounded away from $k$ by $\epsilon'$. So the question is of interest even for a single group $G$.

Let us mention here a result which connects the two topics directly: 
\begin{theorem}[Shalom \cite{shalom-combinatorica}] \label{shal}Let $\Gamma=\langle S \rangle$, $S=S^{-1}$, $|S|=k$, be a finitely generated group and $G=\widehat{\Gamma}$ its profinite completion. Then the Haar integration is the only $\Gamma$-invariant mean on $G$ if and only if the family $\{Cay(\Gamma/N ; S) ; N \triangleleft \Gamma , |\Gamma/N| < \infty \}$ forms an expanding family.
\end{theorem}

In what follows if $G$ is a finite or compact group and $S=S^{-1}$ a subset of $G$ with $|S|=k$, we will say that $S$ is \emph{$\epsilon$-expanding} if all eigenvalues of $\Delta^S = \sum_{s \in S} s$ acting on $L_0^2(G)$ are at most $k-\epsilon$. Sometimes we simply say \emph{expanding} omitting the $\epsilon$, when we talk about an infinite group $G$ or about an infinite collection of $G$'s with the same $\epsilon$. 

In the case of finite groups, if all eigenvalues of $\Delta^S$ on $L_0^2(G)$ are, in absolute value, either $k$ or at most $2 \sqrt{k-1}$, we say that $S$ is a Ramanujan subset of $G$.

In this case $Cay(G,S)$ is a Ramanujan graph (\cite{LPS, lub94}) and $2\sqrt{k-1}$ is the best possible bound one can hope for (for an infinite family of finite groups) by the well-known Alon-Boppana result (see \cite[Prop. 3.2.7]{sarnak-book}). This notion extends also naturally to subsets $S$ of an infinite compact group $G$. Also here $2\sqrt{k-1}$ is the best possible bound (even for a single such group $G$)  because  $2\sqrt{k-1}$ is the rate of exponential growth of the number of closed paths of length $n$ around a base point in the $k$-valent tree. 

In \S 3 and \S 4, we will describe what is known about expanding sets in finite simple groups (in \S 3) and in compact simple Lie groups (in \S 4). Very little is known about the existence of Ramanujan subsets and we will raise there some questions.

\section{Expansion in finite simple groups}

In this section we are interested in expanding subsets of size $k$ in finite simple groups. Abelian groups cannot give rise to expanders, see \cite{lubotzky-weiss}, so when we say simple, we always mean non-commutative simple groups. We will divide our discussion into three subsections: best, random and worst case generators.

\subsection{Best case generators}

The classification of the finite simple groups can be used to show that every such group is generated by two elements. In our context it is natural to ask if all finite simple groups are uniformly expanding. As mentioned in \S 2, it has been shown at an early stage that for fixed $n_0 \geq 3$, $\{\SL_{n_0}(p)\}_{p \in \mathcal{P}}$ is an expanding family when $p$ runs over the set  $\mathcal{P}$ of all primes, using the generators of the ``mother group" $\SL_n(\Z)$. But what about the family $G_n(p_0):=\SL_n(p_0)$ when this time $p_0$ is fixed and $n$ varies?

In \cite{lubotzky-weiss} it was shown that the family $\{G_n(p_0)\}_{n \geq 2}$ is not expanding with respect to \emph{some} set of generators (see \S 3.3 below ; for $\{G_{n_0}(p)\}_{p \in \mathcal{P}}$ this is still an open problem!). This was deduced there by embedding a finitely generated amenable group as a dense subgroup of $\prod_n \SL_n(p_0)$, something which is impossible in $\prod_p \SL_{n_0}(p)$. It has been suggested that maybe ``bounded rank" (i.e. $n_0$ fixed) groups behave differently regarding expansion than unbounded rank (i.e. $p_0$ fixed and $n \to +\infty$). This still might be the case regarding ``worst case generators". A cumulation of works of Kassabov \cite{kassabov-symmetric,kassabov-universal}, Lubotzky \cite{lubJEMS}, Nikolov \cite{nikolov} (c.f \cite{kassabov-lubotzky-nikolov}) and Breuillard-Green-Tao \cite{BGTsuzuki} gives:

\begin{theorem}\label{kln-thm} There exist $k \in \N$ and $\epsilon>0$ such that every non-abelian finite simple group $G$ has a subset $S=S^{-1}$ of size $k$ such that $Cay(G;S)$ is an $\epsilon$-expander.
\end{theorem}

The breakthrough for the proof of Theorem 3.1 was the paper of Kassabov \cite{kassabov-universal} where he broke the barrier of the bounded rank to show that $\{\SL_{3n}(p) ; p \in \mathcal{P}, n \in \N\}$ is an expanding family. Rather than describing the exact historical development (which can be found in \cite{kassabov-lubotzky-nikolov}) let us give the conceptual explanation. 

In \cite{ershov-jaikin} it is shown that $\Gamma= E_d(\Z\langle x_1,\ldots,x_\ell\rangle)$ has property $(T)$ for every $d \geq 3$ and $\ell \in \N$, where $\Z\langle x_1,\ldots,x_\ell\rangle$ is the free non-commutative ring on $\ell$ free variables and $E_d(R)$, for a ring $R$, is the group of $d \times d$ matrices over $R$ generated by the elementary matrices $\{I+re_{i,j} ; 1 \leq i \neq j \leq d, r \in R\}$. Now, for every prime power $q$ and every $n\in \N$, the matrix ring  $M_n(\F_q)$ is $2$-generated as a ring, i.e. $\Z\langle x_1,x_2\rangle$ can be mapped onto the $M_n(\F_q)$. This implies that $\Gamma=E_3(\Z\langle x_1,x_2\rangle)$ can be mapped onto $E_3(M_n(\F_q)) \simeq \SL_{3n}(\F_q)$ and hence by Margulis's original result, i.e. Proposition \ref{margulisprop} above, $\{\SL_{3n}(\F_q)\}_{n \geq 1}$ is an expanding family.

Let us take the opportunity to observe that this can be used to answer (positively!) a question asked in \cite{lubotzky-weiss} ; It was asked there whether it is possible to have an infinite compact group $K$ containing two finitely generated dense subgroups $A$ and $B$ such that $A$ is amenable and $B$ has property $(T)$. If $K$ is a compact Lie group, this is impossible because the Tits alternative forces $A$, and hence $G$ and $B$ to have a solvable subgroup of finite index ; but a $(T)$ group which is also amenable must be finite. On the other hand it was shown in that paper that the compact group $\prod_{n \geq 3} \SL_n(\F_p)$ does contain a finitely generated amenable dense subgroup. Hence its quotient $K:=\prod_{n \geq 1} \SL_{3n}(\F_p)$ also has such an $A$. But from the previous paragraph we see that $K$ also has a $(T)$ subgroup $B$. Indeed the diagonal image of $\Gamma=E_3(\Z\langle x_1,x_2\rangle)$, has $(T)$ and must be dense in $K$ because it maps onto each of the non-isomorphic quasi-simple groups $\SL_{3n}(\F_p)$ (Goursat lemma).

Now let us move ahead with expanders. An easy lemma shows that if a group is a bounded product of expanding groups it is also expanding (for different $k$ and $\eps$). Nikolov \cite{nikolov} showed that when the rank is large enough every finite simple group of Lie type is a bounded product of the groups treated by Kassabov, thus extending the result for all high rank. But what about lower rank and first of all $\SL_2$?

Let us observe first that one cannot hope for a proof of the kind of Margulis/Kassabov for the groups $\SL_2(q)=\SL_2(\F_q)$. In fact there is no ``mother group" $\Gamma$ with property $(T)$ which is mapped onto $\SL_2(q)$ for infinitely many prime powers $q$. Indeed if such a group $\Gamma$ exists, then by some standard ultraproduct argument (or elementary algebraic geometry, see e.g. \cite{LMS} ) $\Gamma$ has an infinite representation into $\SL_2(F)$ for some algebraically closed field $F$. However this is impossible as every property $(T)$ subgroup of  $\SL_2(F)$ must be finite (see e.g. \cite[Thm 3.4.7]{lub94} for the proof when $\textnormal{char}(F)=0$, but the same argument works in positive characteristic : any action of a $(T)$ group on a Bruhat-Tits tree or hyperbolic space must fix a point, so it must have compact closure in all field completions ; see also  \cite[Ch. 6 Prop 26]{delaharpe-valette}).

So a different argument is needed ; for $p$ prime, it has been deduced from Selberg's theorem ($\lambda_1 \geq \frac{3}{16}$) that $\{\SL_2(p)\}_{p \in \mathcal{P}}$ are expanding, see e.g. \cite{gamburd-thesis, breuillard-parkcity}. Similar reasoning (using Drinfled instead of Selberg) gives a similar result for $\{\SL_2(p_0^\ell)\}$ where $p_0$ is fixed and $\ell \in \N$, see \cite{morgenstern}. But how to handle them altogether? This was done by Lubotzky \cite{lubJEMS} using a very specific construction of Ramanujan graphs (and Ramanujan complexes). That construction in \cite{LSV2}, made $G=\SL_2(p^\ell)$ into $(p+1)-$regular Ramanujan graphs using a set of $p+1$ generators of the following type:
$\{tct^{-1} ; t \in T\}$  where $c$ is a specific element in $G=\SL_2(p^\ell)$ and $T$ is a non-split torus in $H=\SL_2(p)$. Now by Selberg as above, $H$ is expanding with respect to $2$ generators (and their inverses), say 

$$a=\left(
  \begin{array}{cc}
   1 & 1 \\
    0 & 1\\
  \end{array}
\right), \textnormal{  and } b=\left(
  \begin{array}{cc}
   1 & 0 \\
    1 & 1\\
  \end{array}
\right)
$$
and $G$ with respect to one $H$-conjugate orbit of $c \in G$. From this one deduces that $G$ is an expander with respect to $\{a^{\pm 1},b^{\pm 1},c^{\pm 1}\}$ (see \cite{lubJEMS}) for details and \cite{kassabov-lubotzky-nikolov}  for an exposition. In fact, it is also shown there that one can use the more general Ramanujan complexes constructed in \cite{LSV1, LSV2} to give an alternative proof to Kassabov's result for $\SL_n(q)$, all $n$ all $q$ simultaneously.

Anyway, once we have also $\SL_2(q)$ at our disposal, all finite simple groups are bounded products of $\SL_n(q)$ (all $n$ all $q$) except for two families that still need a special treatment for the expanding problem and to prove Theorem \ref{kln-thm}.

One family is the family of  Suzuki groups ; these finite simple groups (which are characterized by the fact that they are the only finite simple groups whose order is not divisible by $3$), see \cite{glauberman}, do not contain copies of $PSL_2(\F_q)$ and resist all the above methods. They have been eventually resolved by Breuillard-Green-Tao \cite{BGTsuzuki} by random methods, so we postponed their treatment to \S 3.2.

Last but not least is the most important family of finite simple groups $Alt(n)$. They do contain copies of groups of Lie type, but one can show that they are \emph{not} bounded products of such. So a new idea was needed here ; what Kassabov \cite{kassabov-symmetric} did is to consider first $n$'s of the form $n=d^6$ when $d=2^k-1$ for some $k \in \N$. The fact that $\SL_3(\Z\langle x_1,x_2 \rangle)$ has property $(T)$ implies that the direct product $\SL_{3k}(\F_2)^{d^5}$ is an expanding family and he embedded this group into $Alt(n)$ is $6$ different ways. The product of these $6$ copies is still not the full $Alt(n)$, but (borrowing an idea from Roichman \cite{roichman}) he treated separately representations of $Alt(n)$ corresponding to partitions $\lambda=(\lambda_1 \geq \ldots \geq \lambda_\ell)$ of $n$ with $\lambda_1 \leq n-d^{5/4}$ and all the others. The first were treated by appealing to results on ``normalized character values'' and the second types were treated collectively by giving their sum a combinatorial meaning and showing that the action of a bounded product of six copies of that model mixes in few steps. The reader is referred to the full details of this ingenious proof in \cite{kassabov-symmetric} or to the exposition in \cite{kassabov-lubotzky-nikolov}.

All in all the knowledge on ``best case" expansion in finite simple groups is in a pretty good shape, as Theorem \ref{kln-thm} shows, and certainly better than what we will describe in the other five cases. Still some natural problems arise here: 

\bigskip

\noindent {\bf Problems 3.1} (a) Is there a discrete group with property $(T)$ which is mapped onto all finite simple groups of large rank? or on all $Alt(n)$? By the computer assisted recent breakthrough in \cite{KNO} we now know that the group of (order preserving) automorphisms of the free group $F_5$ has property $(T)$. This group is known to surject onto $Alt(n)$ for infinitely many $n$'s, see \cite{gilman}. See also \cite{lubotzky-autFn}.

(b) The proof of Theorem \ref{kln-thm} described above gives an explicit set of generators in all cases except of the Suzuki groups. It would be of interest to cover also this case.  Also Theorem \ref{kln-thm} gives a certain fixed number $k$ of generators, which is bounded but larger than $2$. One hopes to get a proof with smaller sets of generators (perhaps $2$). This is especially of interest for $Alt(n)$. 

(c) We discussed expanding families, i.e. the eigenvalues are bounded away from $k=|S|$. What about Ramanujan families, i.e. families of groups $G_i$ with $|S_i|=k$ s.t. all non-trivial eigenvalues are bounded by $2\sqrt{k-1}$. As of now, only subfamilies of $\{\SL_2(p^\ell) ; p \textnormal{ prime, }\ell \in \N\}$ are known to have such generators, see \cite{LPS} and \cite{morgenstern}. What about $\SL_3(p)$? $Alt(n)$?  In \cite{parzan-survey} Parzanchevski defines Ramanujan directed graphs. Strangely enough, while it is not known how to turn many finite simply groups into Ramanujan graphs, he manages in \cite{parzan-matrix} to turn them into Ramanujan directed graphs!

\subsection{Random generators}
A well-known result, proved by Dixon in \cite{dixon} for the symmetric groups, Kantor-Lubotzky \cite{kantor-lubotzky} for the classical groups, and Liebeck-Shalev \cite{liebeck-shalev} for the exceptional ones, asserts that for every $k \geq 2$, randomly chosen $k$ elements of a finite simple group $G$ generate $G$. This means that:

$$Prob\big( (x_1,\ldots,x_k) \in G^k ; \langle x_1,\ldots, x_k\rangle = G\big) \longrightarrow_{|G| \to +\infty} 1.$$

The basic question of this section is whether they form expanders, namely is there $\epsilon>0$ such that 
\begin{equation}\label{RandomLimit}Prob\big( Expd(G,k,\epsilon)\big) \longrightarrow_{|G| \to +\infty} 1,\end{equation}
where $$Expd(G,k,\epsilon) := \{(x_1,\ldots,x_k) \in G^k ; Cay(G; S) \textnormal{ is } \eps\textnormal{-expanding for } S=\{x_1^{\pm 1},\ldots,x_k^{\pm 1}\}.$$

This is still widely open. The best result as of now is:

\begin{theorem}[Breuillard-Guralnick-Green-Tao \cite{bggt2}]\label{random-gen} For each $k \geq 2$ and $r\geq 1$, there is $\epsilon>0$ such that $(\ref{RandomLimit})$ holds for all finite simple groups $G$ of rank at most $r$.
\end{theorem}

The rate of convergence in $(\ref{RandomLimit})$ is even polynomial in $|G|^{-1}$. In particular this holds for the groups $G=\PSL_n(q)$ when the rank $n-1$ is bounded and $q$ goes to infinity. It also includes the family of Suzuki groups, thus completing Theorem \ref{kln-thm} by showing the existence of some expanding Cayley graph, see \cite{BGTsuzuki} for this special family, a case which was not covered by the Kassabov-Lubotzky-Nikolov methods. The case of $\PSL_2(p)$ in the above theorem was first established by Bourgain and Gamburd in \cite{bourgain-gamburd}.

The method of proof here, pioneered in \cite{helfgott,bourgain-gamburd}  for the family of groups $\{\PSL_2(p)\}_p$, is based on the classification of approximate subgroups of $G$ (see Theorem \ref{approx-thm} below), an important statistical lemma in arithmetic combinatorics, the so-called Balog-Szemer\'edi-Gowers lemma,  and a crucial property of characters of finite simple groups : the smallest degree of their non-trivial irreducible characters is at least $|G|^\delta$, where $\delta>0$ depends only on the rank of $G$.

This property, going back to Frobenius for $\PSL_2(p)$, was established in full generality in a classic paper by Landazuri-Seitz \cite{landazuri-seitz}. It was used  by Sarnak-Xue \cite{sarnak-xue} and in Gamburd's thesis \cite{gamburd-thesis} in the closely related context of spectral gap estimates for the Laplacian on hyperbolic surfaces. It also plays an important role in various combinatorial questions regarding finite groups. It was coined \emph{quasi-randomness} by T. Gowers \cite{gowers}. In particular it implies what is now called \emph{Gowers' trick}, namely the fact that given any finite subset $A$ of a finite simple group, we have $AAA=G$, that is every element of $G$ can be written as the product of three elements from $A$, provided $|A|>|G|^{1-\delta}$, where $\delta>0$ depends only on the rank of $G$. See \cite{babai-nikolov-pyber, breuillard-msri} for proofs of this fact.

The approximate groups mentioned above are by definition subsets $A$ of $G$ such that $AA$ can be covered by a bounded number of translates of $A$. This bound, say $K$, determines the quality of the $K$-approximate subgroup. With this definition $1$-approximate subgroups of $G$ are simply genuine subgroups of $G$. The classification of subgroups of finite simple groups is a vast subject, which of course is part of the Classification of Finite Simple Groups (CFSG). Starting with Jordan's 19th century theorem that finite subgroups of $\GL_n(\C)$ are bounded-by-abelian \cite{jordan,breuillard-jordan} and Dickson's early 20th century classification of subgroups of $\PSL_2(q)$ \cite{dickson-book}, it climaxes with the Larsen-Pink theorem \cite{larsen-pink}, which gives a CFSG-free classification of subgroups of finite linear groups, saying in essence that they are close to being given by the $\F_q$-points of some algebraic group. Regarding approximate groups the main result is as follows:

\begin{theorem}[Classification of approximate subgroups]\label{approx-thm} Let $G$ be a finite simple group, $A \subset G$ a generating subset and $K\geq 1$. If $AA \subset XA$  for some $X \subset G$ with $|X|\leq K$, then either $|A|\leq CK^C$ or $|G|/|A| \leq CK^C$, where $C$ depends only on the rank of $G$. Moreover for all generating subsets  $A \subset G$,
$$|AAA| \geq \min\{|G|,|A|^{1+\delta}\}$$ for some $\delta>0$ depending only on the rank of $G$.
\end{theorem}

For $\PSL_2(p)$ and more generally in rank $1$ the above result  can be established by elementary methods based on the sum-product theorem \`a la Bourgain-Katz-Tao \cite{bourgain-katz-tao}. This was proved by Helfgott for $\PSL_2(p)$ \cite{helfgott} and generalized by Dinai \cite{dinai} to $\PSL_2(q)$ for all $q$. In high rank new ideas were required and although some mileage had been achieved by Helfgott \cite{helfgott-sl3} and Gill-Helfgott \cite{gill-helfgott} for $\PSL_3(p)$ and $\PSL_n(p)$ the solution came after Hrushovski \cite{hrush} proved a very general qualitative version of the above theorem based on a model-theoretic generalization of the Larsen-Pink theorem and ideas from geometric stability theory. The result in the form above was finally proved independently in \cite{pyber-szabo} by Pyber-Szab\'o (all groups) and in \cite{bgt1} by Breuillard-Green-Tao (who had initially only announced it for Chevalley groups) without using any model theory but rather more down-to-earth algebraic geometry in positive characteristic. See \cite{breuillard-standrews} for an exposition.

We now briefly explain the link between Theorems \ref{approx-thm} and \ref{random-gen} following the strategy first developed in \cite{bourgain-gamburd}. We refer the reader to the expository paper \cite{breuillard-standrews} and to the book \cite{tao-expansion} for further details. In order to get a spectral gap it is enough to show that the probability that the simple random walk on the Cayley graph returns to the identity in $O(\log |G|)$ steps is close to $1/|G|$. The Cayley graph is assumed to have large girth (see below Problem 3.2(c)) so during the first $c\log |G|$ steps ($c$ a small constant), the random walk evolves on a tree and we understand it very well (via Kesten's \cite{kesten} theorem in particular). The main point is then to establish further decay of this return probability between $c\log |G|$ steps and $C \log |G|$ steps ($C$ is a large constant). Here the main tool, the \emph{$\ell^2$-flattening lemma} of Bourgain-Gamburd, is a consequence of the celebrated additive-combinatorial  Balog-Szemer\'edi-Gowers lemma (see \cite{tao-vu}). It asserts that this decay does take place at a lower than exponential rate provided the random walk does not accumulate on an approximate subgroup of $G$. Theorem \ref{approx-thm} then kicks in and allows to reduce the proof to showing that the random walk does not accumulate on subgroups of $G$. This last step, which is in fact the main part of \cite{bggt2} is straightforward in rank $1$ but requires several new ideas in high rank, in particular the existence of so-called \emph{strongly dense} free subgroups of simple algebraic groups in positive characteristic, proved in \cite{bggt1} for this purpose.

Finally we mention that the above method also produces expander Cayley graphs for finite groups that are no longer simple, but appear naturally as congruence quotients of Zariski dense subgroups of arithmetic groups (the so-called \emph{thin groups}, which may have infinite co-volume), e.g. $\SL_d(\Z/n\Z)$ when $n$ is no longer assumed to be prime. This is the subject of \emph{super-strong approximation}, for which we refer the reader to the works of Bourgain, Varj\'u and Salehi-Golsefidy \cite{varju, salehi-varju, bourgain-varju}, and its many applications in particular to the affine sieve \cite{bourgain-gamburd-sarnak,salehi-sarnak} (see also surveys \cite{salehi, breuillard-standrews}).

\bigskip

\noindent {\bf Problems 3.2} (a) The corresponding problem for the family of alternating groups $Alt(n)$ with $n$ growing to infinity is still wide open. Is there $\epsilon>0$ and $k\geq 2$ such that $(\ref{RandomLimit})$ holds when $G=Alt(n)$, or when $G=Sym(n)$ the full symmetric group?  However looking instead at the random Schreier graphs $Sch(X_{n,r};S)$ of $Sym(n)$ obtained by looking at the action on the set $X_{n,r}$ of $r$-tuples of $n$ elements, it is well-known that $(\ref{RandomLimit})$ holds for some $\epsilon=\epsilon(r)>0$ when $r$ is fixed (e.g. see \cite{lub94, friedman-joux-roichman}), while one would need $r=n$ to get the full Cayley graph of $Sym(n)$. Nevertheless a conjecture of Kozma and Puder  (\cite[Conj. 1.8]{parzan-puder})  asserts that for every generating set $S$ the spectral gap of the Cayley graph $Cay(Sym(n);S)$ ought to be entirely governed by that of the Schreier graph $Sch(X_{n,4};S)$ with $r=4$. This conjecture, if true, would imply that random Cayley graphs of $Sym(n)$ are expanders.

(b) Being an expander implies that the diameter of the Cayley graph is logarithmic in the size of the graph, see \cite{lub94}. However when the rank of the finite simple groups goes to infinity, such as for the family $Alt(n)$, we do not even know whether or not the diameter of a random $k$-regular Cayley graph can be bounded logarithmically in the size of $G$. In the case of $Alt(n)$ however poly-logarithmic bounds have been established (see \cite{babai-hayes,schlage-puchta,helfgott-seress-zuk}).

(c) Girth lower bounds are also relevant to the problem. For finite simple groups of bounded rank it is known that a random $k$-regular graph ($k\geq 2$) has girth at least $c\log |G|$, where $c>0$ depends only on the rank. In other words the group's presentation has no relation of length $<c\log |G|$, see \cite{gamburd-hoory}. As pointed above this was used in the proof Theorem \ref{random-gen}. However logarithmic girth lower bounds when the rank of the groups is allowed to go to infinity are still completely open even for $Alt(n)$.

(d) What about Ramanujan graphs? Numerical evidence \cite{lafferty-rockmore} hints that random Cayley graphs of $\PSL_2(p)$ may not to be Ramanujan. However it is plausible that they are in fact \emph{almost Ramanujan}, in the sense that for each $\epsilon>0$ with very high probability as $p \to +\infty$ all non-trivial eigenvalues are bounded by $2\sqrt{k-1} + \epsilon$. See \cite{rivin-sardari} where an upper bound on the number of exceptional eigenvalues is established and numerics are given. The same could be said of the family of alternating groups $Alt(n)$ (and perhaps even of the full family of all finite simple groups). Partial evidence in this direction is provided by Friedman's proof of Alon's conjecture \cite{friedman} that the Schreier graphs of $Alt(n)$ acting on $n$ elements are almost Ramanujan (see also \cite{puder-alon, bordenave-alon, bordenave-collins}).




\subsection{Worse case generators}

The family of finite simple groups $Alt(n)$ was shown (see \S 3.1) to be a family of expanders with respect to some choice of generators, but it is not with respect to others: e.g. take $\tau=(1,2,3)$ and $\sigma=(1,2,\ldots,n)$ if $n$ is odd and $\sigma=(2,...,n)$ if $n$ is even. Then $Cay(Alt(n) ; \{\tau^{\pm 1}, \sigma^{\pm 1}\})$ are not expanders (see \cite{lub94}).

A similar kind of argumentation can be performed for every family of finite simple groups $\{G_i\}_{i \in I}$ of Lie type with unbounded Lie rank. In \cite{somlai} it was shown that for each such family there is a generating set $S_i$ of $G_i$ of size at most $10$ such that the sequence of graphs $\{Cay(G_i,S_i)\}_{i \in I}$ is not expanding. 

By contrast we have the  following conjecture (see \cite{breuillard-icm}):

\begin{conjecture} Given $r\in \N$ and $k \in \N$, there exists an $\eps=\eps(r,k)>0$ such that for every finite simple group of rank $\leq r$ and every set of generators $S$ of size $|S|\leq k$, $Cay(G,S)$ is an $\eps$-expander.
\end{conjecture}

Some evidence towards this conjecture is provided by the following result:

\begin{theorem}[\cite{breuillard-gamburd}] There exists a set of primes $\mathcal{P}_1$ of density $1$ among all primes, satisfying the following: there exists $\eps>0$ such that if $p \in \mathcal{P}_1$ and $x,y$ are two generators of $\SL_2(p)$, then for $S=\{x^{\pm 1}, y^{\pm 1}\}$ $Cay(\SL_2(p), S)$ is an $\eps$-expander.
\end{theorem}

The proof uses the uniform Tits alternative proved in \cite{breuillard-zimmer} as well as the same Bourgain-Gamburd method used in the proof of Theorem \ref{random-gen} above. The uniform Tits alternative in combination with the effective Nullstellensatz is used to show that for most primes $p$, the probability of return to the identity (or to any proper subgroup) of the simple random walk on $\SL_2(p)$ after $n=\log p$ steps is at most $p^{-c}$ for some fixed $c>0$ independent of the generating set. This in turn implies a spectral gap via the Bourgain-Gamburd method and Theorem \ref{approx-thm}.

\section{Expansion in compact simple Lie groups}

In this section we are interested in expanding subsets of size $k$ in compact simple Lie groups. Here again, we will divide our discussion into three subsections: best, random and worst case (topological) generators.

\subsection{Best case generators} Here the question is to find a topological generating set with spectral gap : 

\begin{theorem}[Margulis, Sullivan, Drinfeld] Every simple compact Lie group contains a finite topological generating set with spectral gap.
\end{theorem}

Every simple compact Lie group $G$ not locally isomorphic to $\SO(3)$ contains a countable dense subgroup with Kazhdan's property $(T)$. Indeed one can find an irreducible high rank arithmetic lattice in a product $G \times H$ and project it to $G$,  see \cite[III.5.7]{margulis-book}. Any finite generating subset $S$ of this countable $(T)$ group $\Gamma$ will provide an example of topological generating set of $G$ with a spectral gap (in particular the conditions of Theorems \ref{rosen} and \ref{shal} will hold). These observations were made by Margulis \cite{margulis80} and Sullivan \cite{sullivan81}.

However the case of $\SO(3)$ (and its double cover $\SU(2)$) is exceptional : it does not contain any countable infinite group with property $(T)$ (see \cite[Ch. 6 Prop 26]{delaharpe-valette}). So it seems much harder to find a topological generating set with spectral gap. Nevertheless this was achieved by Drinfeld shortly after the work of Margulis and Sullivan  in a one-page paper \cite{drinfeld}. Kazhdan's original proof that lattices in high rank simple Lie groups have property $(T)$ is representation theoretic by nature. It uses heavily the fact that the discrete group is a lattice, so that one can induce unitary representations from the lattice to the ambient Lie group and thus reduce the problem to a good understanding of the representation theory of the Lie group. Drinfeld's idea is similar : the countable dense subgroup of $G$ he uses arises from the group of invertible elements in a quaternion algebra defined over $\Q$, which ramifies at the real place (so that the associated Lie group is locally isomorphic to $\SO(3)$). But the tools to establish the spectral gap are much more sophisticated : namely the Jacquet-Langlands correspondence is used to reduce the question to spectral gap estimates for Hecke operators associated to irreducible $\PGL_2(\Q_p)$ representations arising from automorphic representations on the adelic space $L^2(\PGL_2(\A)/\PGL_2(\Q))$. These estimates follow either from the work of Deligne on the Ramanujan-Peterson conjectures \cite{deligne} or from earlier estimates due to Rankin \cite{rankin}. We refer to Drinfeld's original paper and to the book \cite{lub94} for the details of this argument.

These methods produce some specific (topological) generating sets arising from generators of a lattice in a bigger group. One can be very explicit and write down concrete matrices for the generators. See \cite{lub94,sarnak-book, cv-bourbaki}. Lubotzky-Phillips-Sarnak \cite{LPS2} pushed this to produce many examples of families of topological generators of $\SO(3)$ with optimal spectral gap (i.e. Ramanujan). For example the set $S$ consisting of the three rotations of angle $\arccos(-\frac{3}{5})$ around the coordinate axes of $\R^3$ and their inverses provides such an example\footnote{These generators are called now $V$-gates in the quantum computing literature, see \S 5.}.

In these examples the quality of the gap deteriorates as the dimension tends to infinity. In \cite[\S 2.4]{sarnak-book}, Sarnak gives an inductive construction starting with a set $S$ of size $k$ in $\SO(n)$ with spectral gap at least $\epsilon$, which produces a new set $S'$ in $\SO(n+1)$ of size $2k$ with spectral gap at least $\epsilon/2k$. However  the following is still open : 

\noindent {\bf Problem 4.1} Does there exists $\eps>0$ and $k>0$ and for each $n$ a $k$-tuple of topological generators of $\SO(n)$ with all eigenvalues $<k-\eps$  (i.e. a spectral gap which is uniform in $n$)?

\subsection{Random case} Here the situation is wide open. Sarnak \cite[p. 58]{sarnak-book} asks the question whether for a generic pair of rotations $a,b$ in $\SO(3)$ the corresponding set $S=\{a,b,a^{-1},b^{-1}\}$ has a spectral gap. This question is still open, should generic be taken either in the sense of Lebesgue measure or in the sense of Baire category. 

Regarding the latter an interesting observation was made in \cite{LPS1} : if $G$ is a compact simple Lie group, then for a Baire generic family of  (topological) generating sets $S$ of size $k$, generating a free subgroup of $G$, the Laplace operator $\Delta^S$ on $L^2_0(G)$ has infinitely many exceptional eigenvalues, i.e. eigenvalues of size $>2\sqrt{k-1}$ ; see below at the end of Section 4.

It is also worth mentioning that  generically a $k$-tuple of elements in $G$ generates a free subgroup. This is true both in the sense of Baire category and in the sense of measure (see e.g. \cite{gartside-knight, breuillard-gelander, aoun}).

Another interesting observation in the random (w.r.t Lebesgue) situation was made by D. Fisher \cite{fisher}. He observed that the property of a $k$-tuple $(a_1,...,a_k)$ in $G^k$ to have some positive spectral gap (i.e. $S=\{a_1^{\pm 1}, \ldots, a_k^{\pm 1}\}$ has a spectral gap) is invariant under the group of automorphisms of the free group $F_k$. Indeed if $S$ has a spectral gap so does any other generating subset of the group $\langle S \rangle$ generated by $S$. Now the action of $Aut(F_k)$ on $G^k$ is known to be ergodic when $k\ge 3$ by a result of Goldman \cite{goldman} (case of $\SU(2)$) and Gelander \cite{gelander} (general $G$) ; see also \cite{lubotzky-autFn} for this and general background on $Aut(F_k)$ and its actions. Hence we have the following zero-one law:

\begin{theorem}[Fisher \cite{fisher}]\label{fisher-thm} Let $G $ be a simple compact Lie group.  If $k \ge 3$, then either Lebesgue almost all $k$ tuples have a spectral gap, or Lebesgue almost no $k$-tuple has a spectral gap.
\end{theorem}

In the next paragraph we discuss the new method introduced by Bourgain and Gamburd to establish a spectral gap. 

\subsection{Worse case} Even worse : we do not know even a single example of a finite topological generating set of $\SO(3)$ (or any compact simple Lie group) without a spectral gap.  Nevertheless a breakthrough took place a decade ago when Bourgain and Gamburd produced many more examples of topological generators with spectral gap. They showed in \cite{bourgain-gamburd-SU2} that every topological generating set all of whose matrix entries are algebraic numbers (i.e. in $\SU(2,\overline{\Q})$) have a spectral gap. This has now been generalized, first to $\SU(d)$ by Bourgain-Gamburd themselves \cite{bourgain-gamburd-SUd} then to arbitrary compact simple Lie groups by Benoist and Saxc\'e \cite{benoist-saxce} : 

\begin{theorem}[Bourgain-Gamburd, Benoist-Saxc\'e] \label{bg-bs} Let $G$ be a compact simple Lie group and $S \subset G$ be a finite symmetric subset generating a dense subgroup and such that the $\tr(Ad(s))$ is an algebraic number for every $s\in S$. Then $S$ has a spectral gap.
\end{theorem}

Note that the best case examples mentioned above and produced by Margulis, Sullivan and Drinfeld have algebraic entries (property $(T)$ groups have algebraic trace field by rigidity) and so they fall in the class of subsets handled by the above theorem. The converse however is not true: in Theorem \ref{bg-bs} the subgroups generated by $S$ are usually not lattices in any Lie group and although they can be made discrete under the usual geometric embedding looking at the different places of the trace field, they will only be \emph{thin subgroups} there, i.e. Zariski-dense of infinite co-volume.

The proof of Theorem \ref{bg-bs} is inspired from the above-mentioned method Bourgain and Gamburd had first pioneered for the family of finite simple groups $\PSL_2(p)$, but it is much more involved. It still contains a significant combinatorial input in that instead of the growth properties of triple products of finite subsets as in Theorem \ref{approx-thm} above, one needs to consider the growth under triple products of $\delta$-separated sets and thus consider \emph{discretized approximate groups}. The cardinality of a finite set $A$ is replaced by the $\delta$-discretized cardinality $\mathcal{N}_\delta(A)$, which is the minimum number of balls of radius $\delta$ needed to cover $A$. This setting was explored by Bourgain in his proof that there is no non-trivial Borel subring of the reals with positive Hausdorff dimension \cite{bourgain-discretized} culminating with \emph{Bourgain's discretized sum-product theorem}, and later by Bourgain-Gamburd in \cite{bourgain-gamburd-SU2}. In his thesis Saxc\'e was able to prove the suitable analogue of Theorem  \ref{approx-thm} in the context of discretized sets in compact Lie groups. 

\begin{theorem}[Saxc\'e's product theorem \cite{saxceproduct}] \label{product-compact} Let $G$ be a simple compact Lie group and $\delta>0$. For every $\kappa>0$ and $\sigma>0$ there is $\epsilon>0$ such that for every set $A\subset G$ which is $(i)$ of intermediate size (i.e. $\mathcal{N}_\delta(A) =\delta^{-\alpha}$ for $\alpha \in [\sigma, \dim G - \sigma]$, $(ii)$ $\kappa$-non-concentrated (i.e. $\mathcal{N}_\rho(A) \geq \delta^\eps \rho^{-\kappa}$ for all $\rho\geq \delta$) and $(iii)$  $\epsilon$-away from subgroups (i.e. for every proper closed subgroup $H$ of $G$ there is $a \in A$ with $d(a,H)\geq \delta^\eps$), we have:
$$\mathcal{N}_\delta(AAA) \geq \mathcal{N}_\delta(A)^{1+\epsilon}.$$
\end{theorem}

An interesting consequence of this theorem (proved in \cite{saxce-duke}) is that sets of positive Hausdorff dimension have bounded covering number : namely given $\sigma>0$ there is $p\in \N$ such that for any Borel subset $A$ of $G$ with Hausdorff dimension at least $\sigma$, $A^p=G$. 

As in the case of finite groups of Lie type the spectral gap in Theorem \ref{bg-bs} is established by showing the fast equidistribution of the simple random walk on $G$ induced by $S$. A similar combinatorial argument based on the Balog-Szemer\'edi-Gowers lemma shows that the fast equidistribution must take place unless the walk is stuck in a $\delta$-discretized approximate group, and hence (applying  Theorem \ref{product-compact}) in the neighborhood of some closed subgroup. It only remains to show that the random walk cannot spend too much time close to any subgroup. This is where the algebraicity assumption comes in. In fact, as shown in \cite{benoist-saxce}, it is enough to know that the probability of being exponentially close to a closed subgroup is uniformly exponentially small. 

\begin{definition}(Weak diophantine property) A finite symmetric subset $S$ of a compact simple Lie group $G$ (with bi-invariant metric $d(\cdot,\cdot)$)  is said to be \emph{weakly diophantine} if there are $c_1,c_2>0$ such that for all large enough $n$ and every proper closed subgroup $H \leq G$ we have 
$$Prob_{\{w  ;|w|=n\}} \big( d(w,H) < e^{-c_1n} \big) \leq e^{-c_2 n}$$
where the probability is taken uniformly over the $k^n$ words $w$ of length $|w|=n$ in the alphabet $S$. 
\end{definition}

Note that the presence of a spectral gap gives a rate of equidistribution of the random walk. In particular,  it easily implies the weak diophantine property. But we now have :

\begin{theorem}[see \cite{benoist-saxce}] For a finite symmetric subset $S$ the weak diophantine property and the existence of a spectral gap are equivalent.
\end{theorem}

If $S$ has  algebraic entries (or if the trace field, generated by $\tr(Ad(s))$ is a number field), then it is well-known and easy to show that it satisfies a ``strong diophantine property", namely there is $c_1>0$ such that $d(w,1)>e^{-c_1|w|}$ for every word $w$ in $S$ not equal to the identity. Benoist and Saxc\'e verify that it must also satisfy the weak diophantine property.

The uniform Tits alternative yields a weaker version of the weak diophantine property (where $e^{-c_1n}$ is replaced by $e^{-n^{c_1}}$) which holds for \emph{every} topological generating set $S$, see \cite{breuillard-zimmer}.

It has been conjectured in \cite{gamburd-jacobson-sarnak}, \cite[\S 4.2]{gamburd} that the strong diophantine property ought to hold for Lebesgue almost every $S$ (it does not hold for every $S$, for example it fails if $S$ contains a rotation whose angle mod $\pi$ is a Liouville number).

Finally we propose a stronger spectral gap conjecture: 

\begin{conjecture} Let $G$ be a simple compact Lie group and $k \geq 4$. There is $\epsilon>0$ such that for every symmetric set $S$ of size $k$ generating a dense subgroup of $G$, $\Delta^S$ has only finitely many eigenvalues outside the interval $[-k+\epsilon, k-\epsilon]$. 
\end{conjecture}

It is easy to see that not all eigenvalues can be contained in a proper sub-interval: for example if the generators are close to the identity in $G$, then there will be many eigenvalues close to the maximal eigenvalue $k$. Partial evidence for this conjecture is supported by the fact that the analogous statement (even without exceptional eigenvalues) does hold, with a uniform $\epsilon$, for the action of $S$ on the regular representation $\ell^2(\langle S \rangle)$ of the abstract group $\langle S \rangle$, as follows from the uniform Tits alternative, see \cite{breuillard-zimmer}.

What about Ramanujan topological generating sets? As mentioned above Lubotzky-Phillips-Sarnak  produced such examples in \cite{LPS}. In \cite[Thm 1.4]{LPS1} however they observed that generic (in the sense of Baire) generators are in general not Ramanujan. In this vein the following is still open :

\bigskip
\noindent {\bf Problem 4.3:} Is being \emph{asymptotically Ramanujan} a  Baire generic property? Namely is there a countable intersection $\Omega$ of dense open subsets of $G^k$ such that for every $\epsilon>0$ and every $k$-tuple of symmetric generators $S \in \Omega$, there are only finitely many eigenvalues of $\Delta^S$ outside the interval $[-2\sqrt{k-1} - \epsilon, 2\sqrt{k-1}+\epsilon]$? The same question can be asked for almost every $k$-tuple of generators in the sense of Lebesgue measure, and in this case the argument for the zero-one law of Theorem \ref{fisher-thm} no longer applies.

\section{Navigation and golden gates}

\subsection{Navigation in finite simple groups} \label{nav-finite}

One of the most important applications of expander graphs is their use for the construction of communication networks. Let us imagine $n$ microprocessors working simultaneously within a super computer. Ideally we would like to have them all connected to each other, but this would require $\Omega(n^2)$ connections, which is not feasible. Expander graphs give a replacement which can be implemented with $O(n)$ connections and with reasonable performances. But this requires also a navigation algorithm, which will find a short path between any two vertices of the graph. It is easy to see that $\epsilon$-expander $k$-regular graphs $X$ have diameter bounded by $C \log_{k-1}(|X|)$, where $C$ depends only on the expansion coefficient $\epsilon$. In \S 3, we showed that there exist $k,\epsilon$ such that every (non-abelian) finite simple group $G$ has a symmetric set of generators $S$ of size at most $k$ such that $Cay(G,S)$ is an $\epsilon$-expander and so there is $C$ such that $\diam(Cay(G,S)) \leq C \log_{k-1}(|G|)$. But the proof that provided these generators did not offer an algorithm to find a path between two given points of length at most $C \log_{k-1}(|G|)$. This is still open : 
\bigskip

{\bf \noindent Problem 5.1.a.}  Find $k \in \N$ and $\epsilon,C>0$ such that every non-abelian finite simple group has a symmetric set of generators of size at most $k$, for which $Cay(G,S)$ is an $\epsilon$-expander and there exists a polynomial time (i.e. polynomial in $\log_{k-1} |G|$) algorithm which expresses any given element in $G$ as a word in $S$ of length at most $C \log_{k-1}(|G|)$.
\bigskip

In \cite{babai-kantor-lubotzky} a set $S$ of size $14$ was presented (for almost all the finite simple groups) for which $\diam(Cay(G,S))=O(\log |G|)$ with an absolute implied constant, even though these Cayley graphs were not uniform expanders. 

The case when $G$ is the alternating group $Alt(n)$ (or the symmetric group $Sym(n)$) is of special interest : we know a set of generators $S$ which would give expanders (see \cite{kassabov-symmetric} or \cite{KNO}) but they come with no navigation. On the other hand  \cite{babai-kantor-lubotzky} gives such a navigation algorithm, but not expanders.

Another case of special interest is the family of groups $\PSL_2(\F_p)$, where $p$ runs overs the primes. For this family $$S = \left \{
\begin{pmatrix}
   1       &  \pm 1  \\
    0       & 1
\end{pmatrix} , 
\begin{pmatrix}
   1      &  0  \\
   \pm 1       & 1
\end{pmatrix}
 \right \}$$
gives rise to expanders, but the best navigation algorithm with this generating set is due to Larsen \cite{larsen-navigation}. However his probabilistic algorithm gives words of length $O((\log p)^2)$ rather than the desired $O(\log p)$. Here is a baby version of this problem.

{\bf \noindent Problem 5.1.b.}  Find an algorithm to express
$\begin{pmatrix}
   1       &  \frac{p+1}{2}  \\
    0       & 1
\end{pmatrix} = 
\begin{pmatrix}
   1      &  1 \\
   0       & 1
\end{pmatrix}^{\frac{p+1}{2}}$ as a word of length $O(\log p)$ using $\begin{pmatrix}
   1       &  \pm 1  \\
    0       & 1
\end{pmatrix}$ and $\begin{pmatrix}
   1       &  0  \\
    \pm 1       & 1
\end{pmatrix}$.

Let us mention that if one allows to add an extra generator, say $t=\begin{pmatrix}
   2       &  0  \\
    0       & \frac{1}{2}
\end{pmatrix}$, then this is easily done. Denote by $u_x$ the unipotent matrix $\begin{pmatrix}
   1       &  x  \\
    0       & 1
\end{pmatrix}$. For $b=1,\ldots,p-1$, write $b=\sum_{i=0}^r a_i 4^i$ with $r \leq \log_4(p)$ and $0\leq a_i \leq 3$. Then, since $tu_xt^{-1}=u_{4x}$, we get $u_b=u_{a_0} t  u_{a_1} t^{-1}  \cdot \ldots \cdot t^r u_{a_r} t^{-r}$. This means that $u_b$ is a word of length $O(r)=O(\log p)$ using only the letters $u_1$ and $t$. A similar trick for the lower unipotent matrices plus the observation that every matrix in $\SL_2(p)$ is a product of at most $4$ upper and lower unipotent matrices shows:

\begin{proposition} Problem 5.1.a. has an affirmative answer for the family of groups $\{PSL_2(p) ; p \textnormal{ prime}\}$. 
\end{proposition}

A lot of efforts have been made to get a good navigation algorithm in $\PSL_2(q)$, with respect to the $p+1$ generators provided by the LPS Ramanujan graph. Sardari \cite{sardari-diameter, sardari} showed that this problem is intimately related to some deep problems in number theory asking for solutions of some diophantine equations.  Some of these problems are NP-complete and some are solved in polynomial time! In \cite{sardari} a probabilistic polynomial time algorithm is given to navigate $G =\PSL(2,q)$ (with respect to the $p+1$ generators of the LPS Ramanujan graphs) finding a path of length at most $ (3 +o(1)) \log |G|$  between any two points. Such a path is at most 3 times longer than the shortest path between the points. On the other hand he shows essentially that  finding the shortest path between any two points is NP-complete!


\subsection{Navigation in simple compact groups}

Let $G$ be a compact group with bi-invariant metric $d$, where our main interest will be compact simple Lie groups with the metric induced by the Riemannian structure.  In this case the analogous question to those discussed in \S \ref{nav-finite} for finite simple groups is of interest even for a single group $G$  and has the following form.

\bigskip

\noindent {\bf Problem 5.2.} Find a finite subset $S$ of $G$ of size $k$, which generates a dense subgroup $\Gamma$ of $G$ and find an algorithm that given $\epsilon>0$ and $g \in G$ provides a word of $w$ of short length in $S$ with $d(w,g)<\epsilon$.

\bigskip

What do we mean by short? Let $\mu_\epsilon$ be the volume of a ball of radius $\epsilon$ in $G$. We normalize the volume so that $\mu(G)=1$. We want to cover $G$ by balls of radius $\epsilon$ around the words of length $\leq \ell$ in $S$. The number $T$ of such words  satisfies $|T| \leq (2k)^\ell$ and so the best we can hope for is $\ell \leq O_k(\log \frac{1}{\mu_\epsilon})$. For simple Lie groups, $\mu_\epsilon \sim c \epsilon^{\dim G}$, so we can hope for $\ell \leq O_{k,G}(\log \frac{1}{\epsilon})$. Ideally we would like also to have an efficient algorithm that, when $\epsilon$ and $g \in G$ are given, will find 
$w \in \Gamma$ which is $\epsilon$-close to $g$ and will express $w$ as a word of length $O(\log(\frac{1}{\epsilon}))$ in the elements of $S\cup S^{-1}$. This problem for the group $\PU(n)$ (especially $\PU(2)$, but also for larger $n$) is of fundamental importance in Quantum Computing. The elements of $\Gamma$ are usually called the ``gates'' and optimal gates are ``golden gates'' (see \cite{nielsen-chuang},  \cite{parzan-sarnak} and the references therein for more on this). We will not go into this direction here, but just mention that the classical Solovay-Kitaev algorithm works for general ``gates'' (i.e. a subset $S$ as before) but gives a word $w$ which is of polylogarithmic  length $\log(\frac{1}{\epsilon})^{O(1)}$, while the spectral methods to be discussed briefly below work for special choices of gates but give $w$ of smaller length, and sometimes even with almost optimal implicit constant. We refer the reader to Varj\'u's work \cite{varju-compact} for the best known polylogarithmic estimates for general gates in all compact simple Lie groups.

Problem 5.2 has two parts and each one is non trivial : $(a)$ Given $g \in G$ find $w \in \Gamma$ which is short and $\epsilon$-close to $g$, and $(b)$ Express $w$ as an explicit short word (``circuit'') in terms of $S$. 

The work of Ross and Selinger \cite{ross-selinger} gave essentially a solution to both parts for the case $G=\PU(2)$. They observe that every $g \in G$ can be written as a product of $3$ diagonal matrices and showed how to solve the problem for each diagonal matrix. They used the group $\Gamma=\PU_2(\Z[\sqrt{2}][\frac{1}{\sqrt{2}}])$ (which is the first factor projection of the corresponding arithmetic lattice in $\PU_2(\R) \times \PGL_2(K)$, where $K$ is the degree $2$ extension of the field $\Q_2$ of $2$-adic numbers associated to the prime $\sqrt{2}$), and solved also Part $(b)$ with this group.

The work of Parzanchevski and Sarnak \cite{parzan-sarnak} gives a conceptual explanation for this and a vast generalization. They find a number of groups $\Gamma$ that are suitable to achieve this goal : all the $\Gamma$'s are arithmetic lattices, which appear naturally as lattices in $\PU(2) \times \PGL_2(K)$, when $K$ is a local non-archimedean field (see \cite{lub94} for a thorough explanation of this). The projection of $\Gamma$ to $\PU(2)$ gives the desired dense subgroup. But the more interesting point is that the discrete projection to $\PGL_2(K)$ and the action of $\Gamma$ on its associated Bruhat-Tits tree  gives the navigation algorithm that solves Part $(b)$ of the problem. Some special choices of such $\Gamma$'s gives ``super golden gates'', which are essentially optimal. 
 
The work of Evra and Parzanchevski \cite{evra-parzan} takes the story a step ahead by studying the analogous problem for $\PU(3)$. This time, this is done via arithmetic discrete subgroups $\Gamma$ of $\PU(3) \times \PGL_3(\Q_p)$. Again the projection to $\PU(3)$ gives the desired dense subgroup of $\PU(3)$, while the projection to the other factor gives an action of $\Gamma$ on the Bruhat-Tits building, which enables also to solve the navigation problem (in spite of not being a tree). All this solves Part $(b)$ of the problem, but not Part $(a)$. The reader is referred to the last paper for this emerging beautiful theory and for more open questions.

\bigskip

\noindent \emph{Acknowledgements.} The authors are grateful to Nicolas de Saxc\'e and to Peter Sarnak for useful comments. The first author acknowledges support from the ERC (grant no. 617129). The second author acknowledges support from the  ERC, the NSF and the BSF. 

\bibliographystyle{alpha}

\newcommand{\etalchar}[1]{$^{#1}$}
\def\cprime{$'$}

\end{document}